\documentclass[12pt,a4paper]{article}
\usepackage{epsfig}
\usepackage[numbers]{natbib}
\usepackage{amsmath, amsthm, amssymb}
\usepackage{epsfig}
\usepackage[margin=3cm, top=3.5cm, bottom=3.5cm]{geometry}

\newcommand\blfootnote[1]{%
  \begingroup
  \renewcommand\thefootnote{}\footnote{#1}%
  \addtocounter{footnote}{-1}%
  \endgroup
}

\newtheorem{remark}{Remark}[section]

\numberwithin{equation}{section}
\everymath{\displaystyle}
\allowdisplaybreaks
\title{Half-plane Green`s function  for  anti-plane elastic 1D hexagonal QCs}
\date{}
\author{Tsviatko V. Rangelov\thanks{Institute of Mathematics and Informatics, Bulgarian Academy of Sciences, 1113, Sofia, Bulgaria}
\footnotemark[1]
\and George D. Manolis\thanks{Department of Civil Engineering, Aristotle University,
 GR-54124, Thessaloniki, Greece}
\and Petia S. Dineva\thanks{Institute of Mechanics, Bulgarian Academy of Sciences, 1113, Sofia, Bulgaria}
}

\begin{document}

\maketitle
\blfootnote{Corresponding author: T. Rangelov, rangelov@math.bas.bg}
\begin{abstract}
\noindent
This paper presents an analytical derivation of a frequency-dependent fundamental solution plus a Green's function for the uni-dimensional, hexagonal quasicrystal sheet subjected to elastic waves under anti-plane strain conditions. Furthermore, closed-form solutions for the free-fields developing in this sheet for  propagating  shear waves are also obtained. The analysis presented here provides a foundation for understanding quasicrystal dynamic behavior and for advancing relevant computational methods.
\end{abstract}

\vspace{2pt}

\noindent
{\bf Key words.}
1D hexagonal quasicrystals; elastic waves; Bak`s model; fundamental solutions; Green’s functions; free-field motions
\vspace{2pt}
\noindent

{\bf Math. Subj. Class.}
35J08, 35Q74

\section{Introduction}
\label{sec1}
Quasicrystals (QC) were first discovered in 1982 by \citet{SBGC84} and represent a distinct class of materials with quasiperiodic atomic arrangements that are neither periodic like conventional crystals nor disordered like amorphous solids. This unique structure imparts exceptional physical and mechanical properties, including low friction and adhesion, high hardness (\citet{TSOYIM92}), enhanced wear resistance (\citet{CQD16}), improved corrosion resistance, and low electrical and thermal conductivities, making them attractive for advanced applications in the aerospace, automotive, energy and protective coating industries (\citet{KBDF12}). Structurally, QC are intermetallic solids with long-range quasiperiodic order and no conventional crystal periodicity, classified as one-, two-, and three-dimensional (\citet{JN90}). In one-dimensional (1D) QC, atomic arrangements are quasiperiodic along one direction and periodic in the perpendicular plane (\citet{WYHD97}). This work focuses exclusively on homogeneous, elastic 1D hexagonal QC, which are of particular scientific interest, especially in light of the seminal work of \citet{MBCJB85} that investigated their quasi-periodicity.

 Despite the growing interest of the engineering society on mechanical modelling of QCs, the research field on the  fundamental solutions and Green`s functions for quasicrystal elastodynamics remain scarce. The present work addresses this issue by deriving frequency-dependent solutions for homogeneous, elastic 1D QC under anti-plane wave motion for both a finite and an infinitely-extending sheet. These solutions underpin the development of boundary integral equation methods (\citet{MDRW17}) that would have to rely on the availability of families of fundamental solutions and/or Green’s functions for quasicrystals. These methods offer certain advantages, including reduced dimensionality, efficient handling of complex geometries, solution evaluation at selected points without full-domain discretization, semi-analytical accuracy, natural treatment of infinite boundaries, and precise computation of stress concentrations near defects. Leveraging these benefits requires rigorous theoretical development combined with careful numerical implementation. To this end, the paper is organized as follows: Sect. \ref{sec2} introduces the problem statement, while in Sect. \ref{sec3} we derive the frequency-dependent elastodynamic fundamental solution and the half-plane Green`s function. Time-harmonic SH plane wave solutions for both a finite and an infinitely-extending QC sheet are also derived and discussed in Sect. \ref{sec4}. Finally, some concluding remarks are given in Sect.\ref{sec6}.

\section{Problem Statement}
\label{sec2}
In the Cartesian coordinate system $Ox_1x_2x_3$  we considered anti-plane strain wave motion in an infinitely-extending 1D  hexagonal QC sheet by assuming qasiperiodicity along the  axis $Ox_3$. Thus, the field variables are independent of the  $x_3$ coordinate and we have an anti-plane strain, phonon-phason coupling elasticity problem in the frequency domain. We now define
\begin{itemize}
\item[(i)] a phonon displacement field $u_3=u_3(\mathbf{x},\omega)$ along with the corresponding stresses $\sigma_{3i}(\mathbf{x},\omega)$ and strain $\varepsilon_{3i}(\mathbf{x},\omega)=0.5u_{3,i}$  fields, where $\mathbf{x}=(x_1,x_2 )$ and $\omega$ is the frequency of vibration, and
\item[(ii)] a phason displacement field $w_3=w_3 (\mathbf{x},\omega)$, with the corresponding stress $H_{3i}(\mathbf{x},\omega)$ and strain $\omega_{3i}(\mathbf{x},\omega)=w_{3,i}(\mathbf{x},\omega)$ fields. Note that commas indicate spatial derivatives with respect to indices $i=1,2$.
\end{itemize}
The constitutive QC law  is, see \citet{WP08}:
\begin{equation}
\label{eq1}
\sigma_{3i}=2c_{44}\varepsilon_{3i}+R_3\omega_{3i}, \quad  H_{3i}=2R_3\varepsilon_{3i}+K_2\omega_{3i}.
\end{equation}
In the above $c_{44}$ is the elastic constant of the phonon field, $R_3$ is the coupling coefficient of the phonon-phason filds, $K_2$ is the elastic constant of the phason field.
These material constants satisfy the following conditions for a well posed problem:
\begin{equation}
\label{eq2} c_{44}>0, \quad K_2>0, \quad R_3>0, \quad c_{44}K_2-R_3^2>0.
\end{equation}
The equations of dynamic equilibrium in time domain are
\begin{equation}
\label{eq2a}
\sigma_{3j,j}=\rho\frac{\partial u_3}{\partial t^2}, \quad H_{3j,j}=\rho\frac{\partial w_3}{\partial t^2},
\end{equation}
where $\rho$ is the material density and $t$ is the time variable.  For time time-harmonic excitations and by inserting the constitutive law for the stresses, the equations of motion attain the following form:
\begin{equation}
\label{eq3}
\left\{\begin{array}{l} c_{44}\Delta u_3+R_3\Delta w_3+\rho\omega^2u_3=0,
\\
R_3\Delta u_3+K_2\Delta w_3+\rho\omega^2w_3=0,\end{array}\right.
\end{equation}
where $\Delta=\frac{\partial}{\partial x_1}+\frac{\partial}{\partial x_2}$ is the 2D Laplace operator. Finally, on the surface defined by a normal vector $\mathbf{n}=(n_1,n_2)$, the   phonon and phason tractions are defined as follows $t_3=\sigma_{3i}n_i$,  $G_3=H_{3i}n_i$.

The derivation below of fundamental solution, Green`s function and free-field waves in QCs is influenced by the anti-plane isotropic case. The only new tool is the decoupling of the QCs system of equations \eqref{eq3} and \eqref{eq4} with  a help of matrix $\mathbf{Q}$ of Eq. \eqref{eq72}.
\section{Fundamental Solution  and  Green`s Function for Quasicristals}
\label{sec3}
We derive here both a fundamental solution and a Green` function for dynamic loads in the form of elastic wave sweeps, which are specific to a QC sheet of infinite and of semi-infinite extent, respectively.

\subsection{Fundamental solution derivation}
\label{sec3-1}
 The fundamental solution of Eqs. \eqref{eq5} is a matrix valued function $\mathbf{v}^\ast=\left(\begin{array}{ll} u^\ast_{31} & u^\ast_{32}\\w^\ast_{31} & w^\ast_{32}\end{array}\right)$, that satisfies the following system  of equations:
\begin{equation}
\label{eq4} \left\{\begin{array}{l}
c_{44}\Delta u^\ast_{3j}+R_3\Delta w^\ast_{3j}+\rho\omega^2u^\ast_{3j}
=-\delta(\mathbf{x}-\mathbf{\xi})\delta_{j1},
\\[1pt]
\\R_3\Delta u^\ast_{3j}+K_2\Delta w^\ast_{3j}+\rho\omega^2w^\ast_{3j}
=-\delta(\mathbf{x}-\mathbf{\xi})\delta_{j2}.
\end{array}\right.
\end{equation}
Both field $\mathbf{x}$ and source $\mathbf{\xi}$ points lie on the $Ox_1x_2$ plane, all indices range as $i,j=1,2$, while $\delta(\textbf{x}-\mathbb{\xi})$ is Dirac`s delta function and $\delta_{ji}$ is the Kroneker delta.
The stresses corresponding to the fundamental solution are denoted as $\sigma_{3ij}$ for the phonon field and $H_{3ij}$ for the phason field:
\begin{equation}
\label{eq5}
\sigma_{3ij}=2c_{44}\varepsilon_{3i,j}+R_3w_{3i,j}, \quad H_{3ij}=2R_3\varepsilon_{3i,j}+K_2w_{3i,j},
\end{equation}
with $\varepsilon_{3ij}=\frac{1}{2}u^\ast_{3i,j}$ and  $\omega_{3ij}=w^\ast_{3i,j}$ the corresponding strains.
Equations \eqref{eq4} can be written in matrix form as
\begin{equation}
\label{eq7}
[\textbf{C}\Delta+\mathbf{\Gamma}]\textbf{v}^\ast(\textbf{x},\mathbf{\xi},\omega)=-\delta(\textbf{x}-\mathbf{\xi})\textbf{I}_2,
\end{equation}
where $\mathbf{C}=\left(\begin{array}{ll}c_{44}&R_3\\R_3&K_2\end{array}\right)$ is symmetric and positive definite, $\mathbf{\Gamma}=\left(\begin{array}{ll}\rho\omega^2&0\\0&\rho\omega^2\end{array}\right)=\rho\omega^2I_2$, and  $\mathbf{I}_2$  is the unit $2\times2$ matrix.

In order to solve Eq. \eqref{eq7}, we find its
canonical form by diagonalizing  matrix $\mathbf{C}$, see \citet{Ge61}.
We start with the eigenvalues $a_{1,2}$ of matrix $\mathbf{C}$, which are positive due to constraint Eq. \eqref{eq2}
\begin{equation}
\label{eq71}
a_{1,2}=\frac{1}{2}[\hbox{Tr}\mathbf{C}\pm\sqrt{(\hbox{Tr}\mathbf{C})^2-4\det
\mathbf{C}}].
\end{equation}
We have that $\hbox{Tr}\mathbf{C}=c_{44}+K_2$  is the trace and $\det \mathbf{C}=c_{44}K_2-R_3^2$  is the determinant of matrix
$\mathbf{C}$. The two corresponding
normalized eigenvectors $|\mathbf{q}_i|=1$ of $\mathbf{C}$ are $\mathbf{q}_1=(q_1^1,q_1^2)$, $\mathbf{q}_2=(q_2^1,q_2^2)$
\begin{equation}
\label{eq73}
\left\{\begin{array}{l} \mathbf{q}_1=\frac{1}{\sqrt{(c_{44}-a_1)^2+R^2_3}}\left(R_3,c_{44}-a_1\right),
\\
\mathbf{q}_2=\frac{1}{\sqrt{(K_2-a_2)^2+R^2_3}}\left( K_2-a_2,R_3\right)
\end{array}\right..
\end{equation}
A compact form emerges by introducing angle  $\psi=\arccos\frac{R_3}{\sqrt{(c_{44}-a_1)^2+R_3^2}}$, so that  $q_1^1=q_2^2=\cos\psi,q_1^2=-q_2^1=\sin\psi$
Then the orthogonal matrix $\mathbf{Q}$ can be express as
\begin{equation}
\label{eq72}
\mathbf{Q}=\left(\begin{array}{ll}q^1_1 & q^1_2 \\q^2_1 &
q^2_2\end{array}\right)=\left(\begin{array}{ll}\cos\psi & -\sin\psi \\\sin\psi &
\cos\psi\end{array}\right),
\end{equation}
and transforms the basis to canonical form as
$\mathbf{Q}^{-1}\mathbf{C}\mathbf{Q}=\left(\begin{array}{ll}a_1& 0 \\0 &
a_2\end{array}\right).$

The next step is to define  a $2\times2$ matrix function $\mathbf{h}$ such that  $\mathbf{v}^\ast=\mathbf{Q}\mathbf{h}$, with $\mathbf{h}$ satisfies the equation
\begin{equation}
\label{eq8} [\mathbf{C}\mathbf{Q}\Delta+\mathbf{\Gamma}
\mathbf{Q}]\mathbf{h}=-\delta(\mathbf{x}-\mathbf{\xi})\mathbf{I}_2.
\end{equation}
Since $\mathbf{\Gamma}=\mathbf{Q}^{-1}\mathbf{\Gamma} \mathbf{Q}$, a left-hand side multiplication of Eq. \eqref{eq8} by  $\mathbf{Q}^{-1}=\left(\begin{array}{ll}\cos\psi & \sin\psi \\-\sin\psi &\cos\psi\end{array}\right)$ yelds
\begin{equation}
\label{eq9} [\mathbf{A}\Delta+\mathbf{\Gamma}
]\mathbf{h}=-\mathbf{Q}^{-1}\delta(\mathbf{x}-\mathbf{\xi}),
\end{equation}
 where $A=\left(\begin{array}{ll}a_1&0\\0&a_2\end{array}\right)$.

Equation \eqref{eq9} consists of four partial differential
equations of Helmholtz type
\begin{equation}
\label{eq10}
[a\Delta+\rho\omega^2]p=-\delta(\mathbf{x}-\mathbf{\xi})q,
\end{equation}
with solution (\citet{Vl71}) is
$p(\mathbf{x},\mathbf{\xi},\omega)=\frac{q}{2\pi a}K_0(-ikr)=\frac{iq}{4a}H^{(1)}_0(kr)$,
 where the radial distance between source and receiver is $r=\sqrt{(x_1-\xi_1)^2+(x_2-\xi_2)^2}$, $k=\sqrt{\frac{\rho}{a}}\omega$ is the wave number
and $K_0, H^{(1)}_0$ are the zero order Kelvin and  first kind, zero order Hankel functions, respectively.

Finally  the fundamental solution $\mathbf{v}^\ast=\left(\begin{array}{ll} u^\ast_{31} & u^\ast_{32}\\w^\ast_{31} & w^\ast_{32}\end{array}\right)$ has the form
\begin{equation}
\label{eq12}\begin{array}{l}
v^{\ast}=\mathbf{Q}\mathbf{h}=\left(\begin{array}{ll} \cos\psi&-\sin\psi\\
\sin\psi&\cos\psi \end{array}\right)
\left(\begin{array}{ll}h^u_{31}&h^u_{32 }\\
h^w_{31}&h^w_{32}\end{array}\right)
\\
=\left(\begin{array}{ll}\cos\psi h^u_{31}-\sin\psi h^w_{31} & \cos\psi h^u_{32}-\sin\psi h^w_{32}
\\
\sin\psi h^u_{31}+\cos\psi h^w_{31}&\sin\psi h^u_{32}+\cos\psi h^w_{32}\end{array}\right)
\end{array}
\end{equation}
with
\begin{equation}
\label{eq12-1}
\begin{array}{ll}
&h^u_{31}=\frac{\cos\psi}{2\pi a_1}K_0(-ik_1r), \quad h^u_{32}=\frac{\sin\psi}{2\pi a_1}K_0(-ik_1r),
\\
&h^w_{31}=\frac{-\sin\psi}{2\pi a_2}K_0(-ik_2r), \quad h^w_{32}=\frac{\cos\psi}{2\pi a_2}K_0(-ik_2r).
\end{array}
\end{equation}

The corresponding stresses and tractions are obtained from Eq. \eqref{eq5} as follows:
\begin{equation}
\label{eq13}\left\{\begin{array}{l}
t_{3i}=\sigma_{3ij}n_j, \quad \sigma_{3ij}=c_{44}u^\ast_{3i,j}+R_3w^\ast_{3i,j},
\\
G_{3i}=H_{3ij}n_j(x), \quad H_{3ij}=R_3u^\ast_{3i,j}+K_2w^\ast_{3i,j}.
\end{array}\right.
\end{equation}
The strains $u^\ast_{3i,j}$ and $w^\ast_{3i,j}$ are computed below taking into account that the derivative $K'_0(z)=-K_1 (z)$
\begin{equation}
\label{eq13-1}\begin{array}{l}
\mathbf{v}^{\ast}_{,j}=\mathbf{Q}\mathbf{h}_{,j}=\left(\begin{array}{ll} \cos\psi&-\sin\psi\\
\sin\psi&\cos\psi \end{array}\right)
\left(\begin{array}{ll}h^u_{31,j}&h^u_{32,j}\\
h^w_{31,j}&h^w_{32,j}\end{array}\right)
\\
=\left(\begin{array}{ll}\cos\psi h^u_{31,j}-\sin\psi h^w_{31,j} & \cos\psi h^u_{32,j}-\sin\psi h^w_{32,j}
\\
\sin\psi h^u_{31,j}+\cos\psi h^w_{31,j}&\sin\psi h^u_{32,j}+\cos\psi h^w_{32,j}\end{array}\right)
\end{array}
\end{equation}
with
\begin{equation}
\label{eq13-2}
\begin{array}{ll}
&h^u_{31,j}=-ik_1\frac{\cos\psi}{2\pi a_1}\frac{r_j}{r}K_1(-ik_1r), \quad h^u_{32,j}=-ik_1\frac{\sin\psi}{2\pi a_1}\frac{r_j}{r}K_1(-ik_1r),
\\
&h^w_{31,j}=ik_2\frac{\sin\psi}{2\pi a_2}\frac{r_j}{r}K_1(-ik_2r), \quad h^w_{32,j}=-ik_2\frac{\cos\psi}{2\pi a_2}\frac{r_j}{r}K_1(-ik_2r).
\end{array}
\end{equation}
Thus, the fundamental solution and its corresponding stresses and tractions for the 1D hexagonal QC depend on its material properties, frequency $\omega$ of the applied load and radial distance $r=|\mathbf{x}-\mathbf{\xi}|$   between source and receiver points.

\subsection{Green`s function derivation}
\label{sec3-2}
The mathematical definition of a ‘fundamental solution’ is the solution of a differential equation with a singular Dirac`s function at the right-hand side and just a radiation boundary condition, while the term ‘Green`s function’ is understood as a fundamental solution that satisfies specific boundary conditions. A Green`s function, if it exists, is unique, while a fundamental solution is not - if we add any   solution of a homogeneous differential equation we obtain again a fundamental solution.

We define a half-plane as $\mathbb{R}^2_-=\{x=(x_1,x_2), x_2<0\}$ and seek a solution $g^\ast$ to Eq. \eqref{eq4}  which satisfies traction-free boundary condition on $x_2=0$,
\begin{equation}
\label{eq12-3}
t^{g^\ast}_{3j}|_{x_2=0}=c_{44}g^{\ast u}_{3j,2}+R_3g^{\ast w}_{3j,2}=0, \quad H^{g^\ast}_{3j}|_{x_2=0}=R_3g^{\ast u}_{3j,2}+K_2g^{\ast w}_{3j,2}=0.
\end{equation}
Thus, a correction term will be added to the fundamental solution of \eqref{eq4} in the form
\begin{equation}
\label{eq6-1}
\mathbf{y}=\left(\begin{array}{ll}u_{31}&u_{32}
\\w_{31}&w_{32}\end{array}\right),
\end{equation}
This correction term is the solution of Eq. \eqref{eq3} and the tractions corresponding to the Green`s function $\mathbf{g}^\ast(\mathbf{x}, \mathbf{\xi},\omega)=\mathbf{v}^\ast(\mathbf{x}, \mathbf{\xi},\omega)+\mathbf{y}(\mathbf{x}, \mathbf{\xi},\omega)$  now satisfy Eq. \eqref{eq4}. We use the notation $t^{g^\ast}$ and $G^{g^\ast}$  for the tractions of the phonon and phason fields corresponding to the half-plane Green`s function. Viewed alternatively,  we are modifying the  fundamental solution $\mathbf{v}^\ast(\mathbf{x},\mathbf{\xi},\omega)$ so that its  tractions $t^{g^\ast}_{3j}|_{x_2=0}$   and $G^{g^\ast}_{3j}|_{x_2=0}$  are zero along the half-plane boundary.
Since the fundamental solution $\mathbf{v}^\ast$ is a $2\times2$ matrix valued function, see  Eqs. \eqref{eq12} and \eqref{eq12-1} that depends only on $r=\sqrt{(x_1-\xi_1)^2+(x_2-\xi_2)^2}$ and $\omega$, for  $\mathbf{x},\mathbf{\xi} \in \mathbb{R}^2_-$ we can denote it as $\mathbf{v}^\ast(r,\omega)$. We can  find the Green`s function using the method of images as in the isotropic anti-plane case. For this, denote $\tilde{r}=\sqrt{(x_1-\xi_1)^2+(x_2+\xi_2)^2}$ and let us define the correction term as $\mathbf{y}=\mathbf{v}^\ast(\tilde{r},\omega)$. It is clear that $2\times2$ matrix function $\mathbf{y}$ in Eq.\eqref{eq6-1} satisfy the Eq. \eqref{eq3}.

Since $\frac{\tilde{r}_2}{\tilde{r}}=-\frac{r_2}{r}$ we have that $v^\ast_{,2}|_{x_2=0}=-y_{,2}|_{x_2=0}$ and correspondingly Eq. \eqref{eq12-3} holds.

In sum, both the fundamental solution and the half-plane Green`s function depend on the material properties $c_{44}, R_3, K_2, \rho$ of the 1D hexagonal QC, on the frequency $\omega$ of the applied load and the distance $r$  between source and receiver points, and distance $\tilde{r}$  between source and image  of the receiver point.
\begin{remark}
\label{rem1}\rm
Consider the case $R_3=0$, then the system of Eq. \eqref{eq3} decouples and we have two anti-plane isotropic equations for $u_3$ and for $w_3$. The eigenvalues in Eq. \eqref{eq71}  are $a_1=c_{44}$, $a_2=K_2$ and wave numbers are: $k_1=\sqrt{\frac{\rho}{c_{44}}}$, $k_2=\sqrt{\frac{\rho}{K_2}}$. Since now $\psi=0$, then in Eq. \eqref{eq72}  matrix $\mathbf{Q}=\mathbf{I}_2$, then the fundamental solution $v^\ast$ in  Eq. \eqref{eq12}, using Eq. \eqref{eq12-1} has the form
\begin{equation}
\label{eq12-r}
\mathbf{v}^{\ast}=\left(\begin{array}{ll} u^\ast_{31} & u^\ast_{32}\\w^\ast_{31} & w^\ast_{32}\end{array}\right)=\frac{1}{2\pi}\left(\begin{array}{ll}\frac{1}{c_{44}}K_0(-ik_1r)&0 \\ 0&\frac{1}{K_2}K_0(-ik_2r)\end{array}\right)
\end{equation}
These solutions $u^\ast_{31}$ and $w^\ast_{32}$  coincide with solutions for the isotropic anti-plane case with constants $c_{44}$ ( for the phonon field) and $K_2$ (for the phason field) respectively derived in many works, see for example \citet{Vl71}, \citet{Ko85}, \citet{Do93}.

Concerning the Green`s function $\mathbf{g}^\ast(x,\xi,\omega)$ for the case $R_3=0$ we obtain
\begin{equation}
\label{eq12-g}\begin{array}{l}
\mathbf{g}^\ast(x,\xi,\omega)=\mathbf{v}^\ast(r,\omega)+\mathbf{v}^\ast(\tilde{r},\omega)
\\
=\frac{1}{2\pi}\left(\begin{array}{ll}\frac{1}{c_{44}}(K_0(-ik_1r)+K_0(-ik_1\tilde{r}))&0 \\ 0&\frac{1}{K_2}(K_0(-ik_2r)+K_0(-ik_2\tilde{r}))\end{array}\right)\end{array}
\end{equation}
These Green`s functions coincide with Green`s function for the isotropic anti-plane case with constants $c_{44}$ (for the phonon field) and $K_2$ (for the phason field) respectively shown,  for example in \citet{Ko85}.
\end{remark}

\section{Free-field Waves in Quasicrystals}
\label{sec4}
We study here the passage of incident waves in QC of both infinite and finite extent. Note that this is scalar wave propagation, albeit coupled, with the closest equivalent being transverse waves in a string. In this classical case, there is a forward moving wave, and upon meeting a boundary generates a reflected backward moving wave.
\subsection{Full-plane wave solutions}
\label{sec4-1}
We have a freely propagating harmonic shear (S) wave in the QC and are asking for a solution in the form
\begin{equation}
\label{eq14}
\left(\begin{array}{l} u^{in}_3\\w^{in}_3 \end{array}\right)=A\left(\begin{array}{l} \zeta_1\\ \zeta_2 \end{array}\right)e^{ik(x_1\cos\varphi+x_2\sin\varphi)}
\end{equation}
where $\varphi\in(0,\pi/2)$ is the incident angle in respect to $Ox_2$ axis. We have $A$ the amplitude of this incident wave, $(\zeta_1,\zeta_2)$  the direction of the wave propagation  and  $k$ is the wave number all to be determined. The vector function Eq. \eqref{eq14} satisfies equation of motion  \eqref{eq3}, i.e.,
\begin{equation}
\label{eq14-1}
\left[\left(\begin{array}{ll}c_{44}\Delta &R_3\Delta \\ R_3\Delta & K_2\Delta\end{array}\right)+\rho\omega^2\mathbf{I}_2\right]\left(\begin{array}{l} u^{in}_3\\w^{in}_3 \end{array}\right)=0.
\end{equation}
After reducing the term  $Ae^{ik(x_1\cos\varphi+x_2\sin\varphi)}$, we obtain
\begin{equation}
\label{eq16}
\left[\left(\begin{array}{ll}-c_{44}k^2+\rho\omega^2 &-R_3k^2 \\-R_3k^2 & -K_2k^2+\rho\omega^2\end{array}\right)\right]\left(\begin{array}{l} \zeta_1\\ \zeta_2 \end{array}\right)=0.
\end{equation}
For  a nontrivial solution of  Eqs. \eqref{eq16} to exists, the  determinant of the system must be set  to zero, i.e., $\det[-\mathbf{C}k^2+\rho\omega^2\mathbf{I}_2]=0$. Now define vector $\left(\begin{array}{l} \tau_1\\ \tau_2 \end{array}\right)$ such that $Q\left(\begin{array}{l} \tau_1\\ \tau_2 \end{array}\right)=\left(\begin{array}{l} \zeta_1\\ \zeta_2 \end{array}\right)$, where matrix $\mathbf{Q}$ is given in Eq. \eqref{eq72}. Left-hand side multiplication of Eq. \eqref{eq16} with $\mathbf{Q}^{-1}$ yelds
\begin{equation}
\label{eq17}
\left[\left(\begin{array}{ll}-a_1 &0  \\0& -a_2\end{array}\right)k^2+\rho\omega^2\mathbf{I}_2\right]\left(\begin{array}{l} \tau_1\\ \tau_2 \end{array}\right)=0.
\end{equation}
In the above $a_1>a_2$ are the two positive eigenvalues of matrix $\mathbf{C}$, see  Eq. \eqref{eq71} and  $a_1$, $a_2$   are defined as the effective shear moduli of the 2 phases. Then
\begin{equation}
\label{eq17-1}
\det\left[\left(\begin{array}{ll}-a_1 &0  \\0& -a_2\end{array}\right)k^2+\rho\omega^2\mathbf{I}_2\right]=(-a_1k^2+\rho\omega^2)(-a_2k^2+\rho\omega^2)=0.
\end{equation}
Finally, we  find two orthogonal vectors as solutions of of  Eq. \eqref{eq17}, for  $k_1=\sqrt{\frac{\rho\omega^2}{a_1}}$, $k_2=\sqrt{\frac{\rho\omega^2}{a_2}}$ respectively:
\begin{equation}
\label{eq17-2}
\tau^1=\left(\begin{array}{l} \tau^1_1\\ \tau^1_2 \end{array}\right)=\left(\begin{array}{l} 1\\ 0 \end{array}\right),  \hbox{ and } \tau^2=\left(\begin{array}{l} \tau^2_1\\ \tau^2_2 \end{array}\right)=\left(\begin{array}{l} 0\\ 1 \end{array}\right),
\end{equation}
then
\begin{equation}
\label{eq17-3}
\left(\begin{array}{ll}\zeta^1_1 &\zeta_1^2 \\ \zeta_2^1& \zeta_2^2\end{array}\right)=Q\left(\begin{array}{ll}1 &0\\ 0& 1\end{array}\right)=\left(\begin{array}{ll}\cos\psi &-\sin\psi \\ \sin\psi& \cos\psi\end{array}\right).
\end{equation}
Since in a waveguide defined by a scalar displacement will support a single propagating wave, either in the positive or negative directions, wavenumbers $k_1$ and $k_2$ must correspond to the photon/phason cases or a combination thereof. These propagate in the infinite, 1D hexagonal QC plane with a fixed incident angle $\varphi\in(0,\pi/2)$ with respect to $Ox_1$ axis. Using the definition of $\cos\psi$ and $\sin\psi$  in Eq. \eqref{eq72} it is easy to check that the incident wave solutions have the form:
\begin{equation}
\label{eq18}
\begin{array}{l} \left(\begin{array}{l} u^{in,S_1}_3\\w^{in,S_1}_3 \end{array}\right)=A\left(\begin{array}{l} \cos\psi\\ \sin\psi \end{array}\right)e^{ik_1(x_1\cos\varphi+x_2\sin\varphi)}
\\
\left(\begin{array}{l} u^{in,S_2}_3\\w^{in,S_2}_3 \end{array}\right)=A\left(\begin{array}{l} -\sin\psi\\ \cos\psi \end{array}\right)e^{ik_2(x_1\cos\varphi+x_2\sin\varphi)}
\end{array}
\end{equation}
Note that in contrast with an elastic, homogeneous and isotropic material, the direction vector is a constant and does not depend on the incident angle $\varphi$. In Eq. \eqref{eq18}, we have two modes for the incident S wave, which are orthogonal linear combinations of the phonon ($u_3$) and phason ($w_3$) displacements. Since $a_1$ and $a_2$ are the effective shear moduli of the two phases, the corresponding wave speeds are $c_1=\sqrt{\frac{a_1}{\rho}}$ and $c_2=\sqrt{\frac{a_2}{\rho}}$  where $c_1>c_2$. Thus, we have the propagation of both a fast ($S_1$) and a slow ($S_2$) shear wave due to the phonon-phason coupling effect.

\subsection{Half-plane wave solutions}
\label{sec4-2}
In deriving the half-plane solutions, we follow \citet{Ac73} for the isotropic anti-plane case. The free-field solutions are sum of incident and reflected waves such that their traction on $Ox_2$ is zero.

Using the Eq. \eqref{eq18} there are  two kind of free-field solutions.

With incident $S_1$  wave:
\begin{equation}
\label{eq18-1}
\left(\begin{array}{l} u^{ff,S_1}_3\\w^{ff,S_1}_3 \end{array}\right)=A_1\left(\begin{array}{l} \cos\psi\\ \sin\psi \end{array}\right)\left(e^{ik_1(x_1\cos\varphi+x_2\sin\varphi)}+e^{ik_1(x_1\cos\varphi-x_2\sin\varphi)}\right).
\end{equation}
With incident $S_2$ wave:
\begin{equation}
\label{eq18-2}
\left(\begin{array}{l} u^{ff,S_2}_3\\w^{ff,S_2}_3 \end{array}\right)=A_2\left(\begin{array}{l} -\sin\psi\\ \cos\psi \end{array}\right)\left(e^{ik_2(x_1\cos\varphi+x_2\sin\varphi)}+e^{ik_2(x_1\cos\varphi-x_2\sin\varphi)}\right).
\end{equation}
Here, $A_1$ and $A_2$ are the amplitudes of the incident wave.
In closing, the displacement field in a homogeneous  1D elastic hexagonal QC half-plane subjected to anti-plane wave motion is presented by incident time-harmonic fast $S_1$ and slow $S_2$ waves which generate only reflected waves of the same type.  They depend on the wave frequency,  wave propagation direction, phonon-phason coupling and also on material properties of the QCs.

\section{Conclusions}
\label{sec6}
The aim of this work was to derive frequency-dependent fundamental solutions for the homogeneous elastic 1D QC, plus Green`s functions in the presence of a free surface, all under anti-plane wave motion. These analytical results provide tools for understanding the dynamic behavior of QC and form a foundation for subsequent development of advanced computational methods. In particular, these solutions enable the development of mechanical models based on the boundary integral equation method, an efficient and accurate mesh-reduction technique whose advantages are well documented in the literature. Even on their own, these analytical solutions allow for the treatment of a number of problems involving QC with simple geometries.
\vspace{5pt}

\textbf{Acknowledgement:} The authors wish to acknowledge the Joint BG-GR Research Project of the Bulgarian Academy of Sciences entitled “Physics-informed Neural Networks in Elastodynamics based on Green`s Function Solutions” for the period 2025-2027.

\end{document}